\documentclass[11pt]{article}
\usepackage{amsmath,amsthm,amssymb,hyperref}

\newcommand{\al}{\alpha}

\newcommand{\ep}{\epsilon}
\newcommand{\la}{\lambda}

\newcommand{\om}{\omega}
\newcommand{\bR}{\mathbb{R}}
\newcommand{\bZ}{\mathbb{Z}}

\newcommand{\pt}{\partial_t}
\newcommand{\pa}{\partial}

\newcommand{\beeq}{\begin{equation}}
\newcommand{\eneq}{\end{equation}}

\newtheorem{thm}{Theorem}
\newtheorem{prop}{Proposition}
\newtheorem{coro}[prop]{Corollary}

\newtheorem{lem}{Lemma}

\newenvironment{prf}{\noindent {\bf Proof} }{\endprf\par}
\def \endprf{\hfill  {\vrule height6pt width6pt depth0pt}\medskip}

\numberwithin{equation}{section}

\begin{document}

\title{Local Existence for Nonlinear Wave Equation with Radial Data in $2+1$ Dimensions\thanks{Supported by NSF of
China 10571158}}
\author{ Chengbo Wang and Daoyuan Fang
\thanks {email: CW: wangcbo@yahoo.com.cn, DF: dyf@zju.edu.cn}
\\ Department of Mathematics, Zhejiang University,\\ Hangzhou, 310027, China
} \maketitle

\section{Introduction}\label{4-sec-intro}

Let $\Box_{g}=\pt^2 - g \Delta_x$. In this paper, for initial data
$(u_0, u_1)$ with spherical symmetry and $s>\frac{3}{2}$, we
consider the local well-posed (LWP) or local existence result for
the semilinear (with $g(u)\equiv 1$, SLW) and quasilinear wave
equation (with $g(0)=1$, QLW) \beeq
\label{4-qlw}\left\{\begin{array}{l} \Box_{g(u)}
u = p(u) (\partial_t u)^2 + q(u) (\nabla u)^2 :=N(u, \partial u)\\
u(0,x)=u_0\in H^s,\ \pt u(0,x)=u_1\in H^{s-1}\end{array}\right.\eneq
on $\bR\times \bR^2$. We use $\pa$ to stand for space-time
derivatives, i.e. $\pa=(\pt, \pa_x)$.

For general spatial dimensions $n$, the critical index of Sobolev
space for such problem is $s_c=\frac{n}{2}$ and the counterexamples
for LWP give the lower bound 
$\max(\frac{n}{2}, \frac{n+5}{4})$ (see \cite{FW1} and \cite{Lind93}
for example). The classical theory (see \cite{HuKaMar77} for
example) says that this problem is LWP in $H^s\times H^{s-1}$ for
$s>\frac{n+2}{2}$, which insures that
$\pa u$ is bounded.

For semilinear problems, it's known that one can improve the needed
regularity to $s>\max(\frac{n+1}{2}, \frac{n+5}{4})$ with the help
of Strichartz estimates (see \cite{FW1} for example), and the sharp
results have been proved to be $s>\max(\frac{n}{2}, \frac{n+5}{4})$
(see \cite{Ta99} and references therein).

In the last ten years or so, the analysis of QLW has experienced a
dramatic growth. Following partial results independently obtained
by Bahouri-Chemin \cite{BaCh1}, \cite{BaCh2} and Tataru
\cite{Ta00}, \cite{Ta01}, \cite{Ta02}, and further work of
Klainerman-Rodnianski \cite{KlRo03}, Smith and Tataru largely
completes the local theory for general second order quasilinear 
hyperbolic equations in \cite{SmTa05}. They show that for $n\le
5$, the problem is LWP in $H^s$ for $s>\max(\frac{n+1}{2},
\frac{n+5}{4})$. Moreover, the $L^{\max(2,
\frac{4}{n-1})}L^\infty$ Strichartz estimate holds true for the
corresponding wave operator $\Box_{g(u)}$. For the detailed
historical introduction, see Section 1.2 in \cite{SmTa05} for
example.

Thus, in general, the optimal regularity for LWP for the problem in
two space dimensions is $\frac{7}{4}$. Note 
that the counterexample which gives the lower bound $\frac{7}{4}$
is non-radial.  The main purpose of this paper is to show that for
the radial data, the regularity can be improved further in two
space dimensions by showing that we can get an improved Strichartz
estimates.

In our previous paper \cite{FW}, we get the following radial
improvement of Strichartz estimate for the solution of linear wave
equation $\Box u=0$, \beeq\label{4-StriSLW}\|\pa u \|_{L_{loc,t}^2
L_x^\infty}\lesssim \|\pa u(0)\|_{H^{s-1}}\eneq with $s>
\frac{3}{2}$ and $n=2$. This would naturally yield the radial LWP
in $H^s$ with $s> \frac{3}{2}$ for SLW. We will give the proof in
Section \ref{4-sec-slw}. Moreover, we get a weak stability result
in a more larger class, for any space dimensions. Inspired by the
result for SLW, we intend to prove a similar result in the
quasilinear case, by using the method in \cite{SmTa05}.

Now we begin to state our main results.
\begin{thm}[Radial LWP for SLW]
\label{4-thm-lwpslw} Let $n=2$ and  $s>\frac{3}{2}$. The equation
\eqref{4-qlw} with $g(u)\equiv 1$ is radial local well-posed in $C
H^s$. Precisely, for any radial data $(u_0,u_1)\in H^s\times
H^{s-1}$, there exists a unique radial local in time solution
$u\in C H^s$ such that $\pa u\in C H^{s-1}\cap L^2 L^\infty$, and
the solution map is Lipschitz continuous on bounded sets.
\end{thm}

As we know, the current counterexample to radial LWP shows the
lower bound $s_c=\frac{n}{2}$. And here the positive result
requires the regularity $s>\frac{n+1}{2}$. Thus if we combine it
with the previous positive results $s>\max(\frac{n+5}{4},
\frac{n}{2})$, we know that there is still $\frac{1}{2}$ gap for
$n=2, 3$, and $\frac{1}{4}$ gap for $n=4$. Thus a natural problem
is: {\center \bf What is the optimal regularity $s_o$ for SLW to
be radial LWP?} \\ 
We conjecture that $s_o=\frac{n}{2}$. We still can't prove or
disprove the conjecture now, instead, we utilize the energy estimate
to establish the following weak stability estimate for SLW.

Note that for the equation of type $\Box u= u \nabla u $ in four
space dimensions, Sterbenz \cite{Stbz04p} got a relative results of
global existence with small data, based on the argument of Tataru
\cite{Ta99}. We intend to solve the conjecture by similar method in
the following work.

\begin{thm}[Weak Stability for SLW]
\label{4-thm-unique} Let $s>\frac{n}{2}$ and $s\ge 1$. Consider the
semilinear wave equation \beeq \label{4-gslw}\left\{\begin{array}{l}
\pt^2 u - \Delta
u = \sum_{|\al|=2} q_{\al}(u) (\pa u)^{\al}\\
u(0,x)=u_0,\ \pt u(0,x)=u_1\ .\end{array}\right.\eneq then there
exists at most one solution in the solution class $X=\{u\in C H^s;
\pa u\in C H^{s-1}\cap L^1 L^\infty\}$.  Moreover, if $u$, $v$ are
two solutions of above equation with initial data $(u_0,u_1)$ and
$(v_0,v_1)$, then $$\|\pa(u-v)\|_{L^\infty([0,T], L^2)} \lesssim
\|(u_0-v_0, u_1-v_1)\|_{H^1\times L^2}  \exp(C (T+ \|\pa
(u,v)\|_{L^1 L^\infty}) )\ .$$
\end{thm}

\begin{thm}[Local Existence for QLW]\label{4-thm-localexst}
Let $n=2$ and  $s>\frac{3}{2}$. For each $R > 0$, there exist
constants $T,M> 0$ so that, for any radial initial data $(u_0, u_1)$
which satisfies
$$\|(u_0, u_1)\|_{H^s\times H^{s-1}}\le R\ ,$$
there exists a radial solution $u\in C H^s$ to \eqref{4-qlw} on
$[-T,T] \times \bR^2$ such that \beeq\label{4-FianlEst}\|\pa
u\|_{C_t H_x^{s-1}\cap L^2_t L_x^\infty}\le M\ .\eneq Moreover, for
$1\le r \le s+1$, each $t_0\in [-T,T]$, and any radial data $(v_0,
v_1)$,  the linear equation \beeq\label{4-LW}\left\{\begin{array}{l}
\Box_{g(u)}v=0,\ (t,x)\in [-T,T]\times \bR^2\\
v(t_0)=v_0\in H^r,\ \pt v(t_0)=v_1\in H^{r-1}
\end{array} \right.\eneq
admits a radial solution $v\in C([-T,T],H^r)\cap C^1([-T,T];
H^{r-1})$, and the following estimate holds: \beeq\label{4-energy}
\|v\|_{L^\infty_t H_x^r}+\|\pt v\|_{L^\infty_t H_x^{r-1}} \le C
\|(v_0,v_1)\|_{H^r\times H^{r-1}}.\eneq Additionally, the
following estimate holds, provided $\rho<
r-\frac{1}{2}$,\beeq\label{4-StriQLW}\|\langle D_x \rangle^\rho
v\|_{L^2_t L_x^\infty}\le C \|(v_0,v_1)\|_{H^r\times H^{r-1}}\eneq
and the same estimate holds true with $\langle D_x \rangle^\rho$
replaced by $\langle D_x \rangle^{\rho-1}\pa$.
\end{thm}
As in \cite{SmTa05}, for the proof of Theorem \ref{4-thm-localexst},
we will mainly prove the following dispersive (Strichartz) estimate.

\begin{thm}[Dispersive Estimate]\label{4-thm-disp}
Let $\ep_0 \la\gg 1$ and $\chi_{j,k}$ ($j,k\in \bZ$) be the
``radial" wave packet which will be introduced in Section
\ref{4-sec-RedDisp}, \beeq \label{4-mainest}\|\sum
a_{j,k}\chi_{j,k}\|_{L^2_t L^\infty}\lesssim \ep_0^{-\frac{7}{4}}
(\ln\la)^{\frac{3}{2}}\|a_{j,k}\|_{l_{j,k}^2}\eneq
\end{thm}

We give here some notations which will be used hereafter. Let
$\langle x\rangle=\sqrt{1+x^2}$ and $H(x)$ be the usual Heaviside
function($H(x)=1$ if $x\ge 0$ and $H(x)=0$ else). For a set $E$,
we use $|E|$ to stand for the measure or cardinality of the set
$E$ depending on the context.

This paper will be organized as follows: First, for SLW, we give
the proof of radial LWP result(Theorem \ref{4-thm-lwpslw}) in
Section \ref{4-sec-slw}, moreover, we give the proof of
uniqueness and weak stability result (Theorem \ref{4-thm-unique})
in a more larger class, for any space
dimensions, by utilizing the energy estimates.

Then we turn into  the proof of Theorem \ref{4-thm-localexst}. 
In Section \ref{4-sec-Disp2LocalExst}, we reduce Theorem
\ref{4-thm-localexst} to Theorem \ref{4-thm-disp} and
in Section \ref{4-sec-sketch}, we give the sketch of the proof of
Theorem \ref{4-thm-disp}. We will reduce the proof to the
corresponding estimates for two sets
of indices $(j,k)$ for wave packets separately. 

The first case is that $\chi_{j,k}$ is essentially
$L^\infty$-normed. For this case, we get the overlap estimate of the
``radial" wave packets in Section \ref{4-sec-overlap},
then we have the required dispersive estimate as in Section 10 of
\cite{SmTa05}.

For the remained case, 
it turns out that this is the case when the wave packet evolve
essentially along the light cone and occurs only when $j$ is small.
For this case, we give a $L^\infty$ estimate for single $\chi_{j,k}$
in Section \ref{4-sec-est-ra-wp}, which is sufficient for the proof
of the dispersive estimate, as explained in the end of Section
\ref{4-sec-sketch}.

\section{LWP and Weak Stability for SLW}\label{4-sec-slw}

In this section, we prove the results for SLW. First, we prove the
radial local well posed result Theorem \ref{4-thm-lwpslw}.

{\noindent {\bf Proof of Theorem \ref{4-thm-lwpslw}.} } The
existence of the radial solution follows from the radial Strichartz
estimate. Precisely, for radial $u$ such that $\Box u=0$, we have
$$\|\pa u \|_{L_{loc,t}^2 L_x^\infty}\lesssim \|\pa u(0)\|_{H^{s-1}}$$ with
$s> \frac{3}{2}$. We get the solution by contraction argument as
usual. Let $(u_0, u_1)$ be radial and $$\|u_0\|_{H^s}+\|u_1\|_{
H^{s-1}}\le M\ .$$ Define a complete domain with $C$ large enough
$$B_{T}=\{u\in C([0,T], H^s)\cap C^1 H^{s-1}\ |\ u\ \mathrm{radial},\
\|u\|_{L^\infty H^s}+\|\pa u\|_{ L^\infty H^{s-1}\cap L^{2}L^\infty
}\le C M\}\ ,$$ and for $u\in B_{\ep, T}$, define $\Pi(u)$ be the
solution of the equation
$$\Box \Pi(u)=N(u, \partial u)$$ with prescribed initial data $(u_0, u_1)$.

Thus for $T$ small enough, by Strichartz estimate and energy
estimate,
$$
\begin{array}{lcl}
\|\pa \Pi(u)\|_{L^\infty H^{s-1}\cap
L^{2}L^\infty}  &\le& C_1( M+\|N(u,\partial u)\|_{L^1 H^{s-1}})\\
& \le & C_2 (M+ T^{\frac{1}{2}}\|\pa u\|_{L^\infty
 H^{s-1}\cap L^{2}L^\infty}^2)\\
&\le& C_3( M +T^{\frac{1}{2}} (C M)^2)\le \frac{ C M}{4}
\end{array}
$$
and by noting that $u(t)=u_0+\int_0^t \pt u$ and $s>1$,
$$
\begin{array}{lcl}
\|\Pi(u)(t)\|_{H^s} &\le & \|\pa_x \Pi(u)(t)\|_{H^{s-1}}+\|\Pi(
u)(t)\|_{L^2}\\
& \le & \frac{C M}{4}+M+ \|\pt u\|_{L^1 L^2}\\
& \le & (1+T)\frac{C M}{4}+M\le \frac{C M}{2}\ .\end{array}$$ Thus
$\Pi$ is closed in the ball $B_{T}$, similar argument shows that
$\Pi$ is a contraction map in the ball $B_{T}$. So we get a radial
local solution $u\in  C H^s\cap C^1 H^{s-1}$.

It is easy to see that the radial solution is unique and the
solution map is Lipschitz continuous on $B_{T}$ by the previous
argument. {\hfill {\vrule height6pt width6pt depth0pt}\medskip}

Now we give the proof of Theorem
\ref{4-thm-unique}.

{\noindent {\bf Proof of Theorem \ref{4-thm-unique}.} } Let
$s>\frac{n}{2}$, $$X=\{u\in C([0,T], H^s); \pa u\in C H^{s-1}\cap
L^1 L^\infty\}\ ,$$ and $u$, $v$ in $X$ be two solutions of equation
\eqref{4-gslw} with initial data $(u_0,u_1)$ and $(v_0,v_1)$, then
$\om:=u-v \in C H^s$, and
$$\Box \om=a(u,v)\om (\pa u)^2+q(v) \pa( u, v)\pa
\om :=\tilde{N}$$ with initial data $(\om_0,\om_1)$. Note that
$\om(t)=\om(0)+\int_0^t \pt \om$, then $$\|\om\|_{L^\infty
H^1}\lesssim \|\om\|_{L^\infty L^2}+\|\pa \om\|_{L^\infty
L^2}\lesssim \|\om(0)\|_{L^2}+\|\pa\om\|_{L^\infty L^2}$$ with
$T\lesssim 1$. Thus by Leibnitz rule and Sobolev multiplication law,
we have
$$\begin{array}{lcl}\|\pa\om\|_{L^\infty L^2}
&\lesssim &\|\pa\om(0)\|_{L^2}+\|\tilde{N}\|_{L^1 L^2}\\
&\lesssim &\|\pa\om(0)\|_{L^2}+ \|\om\|_{L^\infty H^1} \|(\pa
u)^2\|_{L^1 H^{s-1}}+\|\pa (u, v)\|_{L^1
L^\infty}\|\pa\om\|_{L^\infty L^2}\\
&\lesssim & \|\pa\om(0)\|_{L^2}+ \|\om\|_{L^\infty H^1} \|\pa
u\|_{L^1 L^\infty}\|\pa u\|_{L^\infty
H^{s-1}}\\
&&+\|\pa (u, v)\|_{L^1 L^\infty}\|\pa\om\|_{L^\infty L^2}
\\
&\lesssim &\|(\om_0,\om_1)\|_{H^1\times L^2}+ \|\pa \om\|_{L^\infty
L^2} \|\pa (u,v)\|_{L^1 L^\infty}\ .
\end{array}$$
So we have the following stability estimate for small enough time
$T\in (0,1]$ (such that $\|\pa (u,v)\|_{L^1 L^\infty}\ll 1$ ),
$$\|\pa\om\|_{L^\infty L^2} \lesssim \|(\om_0,\om_1)\|_{H^1\times L^2}\ .$$
Thus by an induction argument we can get the final estimate
$$\|\pa \om\|_{L^\infty([0,T], L^2)} \lesssim \|(\om_0,\om_1)\|_{H^1\times L^2} \exp(C (T+\|\pa (u,v)\|_{L^1 L^\infty} ))\ .$$
{\hfill {\vrule height6pt width6pt depth0pt}\medskip}

\section{Local Existence for Quasilinear Wave Equation}
\label{4-sec-Disp2LocalExst}

In this section we reduce our main result
Theorem \ref{4-thm-localexst} to the dispersive estimate. 

\subsection{Existence Result for Smooth Initial Data}
First, we show that Theorem \ref{4-thm-localexst} is a consequence
of the following existence result for smooth initial data.
\begin{prop}[Local existence for smooth data]\label{4-prop-localexst}
Let $n=2$ and $s>3/2$. For each $R > 0$, there exist constants $T,M>
0$ so that, for any smooth radial data $(u_0, u_1)$ which satisfies
$\|(u_0, u_1)\|_{H^s\times H^{s-1}}\le R$, there exists a unique
smooth solution $u$ to \eqref{4-qlw} on $[-T,T] \times \bR^2$ such
that
$$\|\pa u\|_{C_t H_x^{s-1}\cap L^2_t L_x^\infty}\le M\ .$$ Moreover,
we have the energy estimate \eqref{4-energy} and Strichartz estimate
\eqref{4-StriQLW} for the solution $v$ of the equation $\Box_{g(u)}
v=0$.
\end{prop}
In fact, for any radial initial data $(u_0, u_1) \in H^s \times
H^{s-1}$ such that $$\|(u_0, u_1)\|_{H^s \times H^{s-1}}\le R\ .$$
Let $(u^k_0, u^k_1)$ be a sequence of smooth data converging to
$(u_0, u_1)$, which also satisfy the same bound. Then the conclusion
of Proposition \ref{4-prop-localexst} applies uniformly to the
corresponding solutions $u^k$. In particular, it follows that the
sequence $\pa u^k$ is bounded in the space $ C H^{s-1} \cap L^2
L^\infty$. Thus there exists a subsequence (also denoted by $u^k$)
which converges weakly to some $u$ in $C H^s\cap C^1 H^{s-1}$. We'll
show below that it's a solution of the equation with data $(u_0,
u_1)$.

Let $\phi_j(x)=\phi(j^{-1} x)$, where $\phi$ is a smooth bump
function with compact support, $\phi=1$ on the unit ball. For any
fixed large $j$, define $u_j^k= \phi_j u^k$, and thus $u_j^k$ are
uniformly bounded in $C H^s\cap C^1 H^{s-1}$. Thus by compactness,
there is a subsequence (also denoted by $u_j^k$) which converges
to some $u_j$ in $C H^{s-}\cap C^1 H^{s-1-}$. However, since
$u_j^k=u^k$ in $B_j$, we have $u_j=u$ in $B_j$.

As a consequence of the fractional Leibnitz rule, the right hand
side term $N(u^k,\partial u^k)$ of the equations for $u^k$ are
uniformly bounded in the space $L^2 H^{s-1}$. Then
\eqref{4-StriQLW} combined with Duhamel's formula show that $\pa
u^k$ is uniformly bounded in $L^2 C^\delta$. Note that
$s>\frac{n}{2}$ and $$\pa u_j^k=\phi_j(x) \pa u^k+j^{-1}(\pa
\phi)(j^{-1} x) u^k\ ,$$ thus we have $\pa u_j^k$ is also
uniformly bounded in $L^2 C^\delta$. Together with the above this
implies that $\pa u_j^k$ converges to $\pa u_j$ in $L^2 L^\infty$.
Thus we get that $\pa u^k$ converges to $\pa u$ in $C
H_{loc}^{s-1-}\cap L^2  C_{loc}^\delta$.

The above information is more than sufficient to allow passage to
the limit in the equation \eqref{4-qlw} and show that $u$ is a
solution in the sense of distributions, yielding the existence part
of Theorem \ref{4-thm-localexst}. The conditions \eqref{4-FianlEst},
\eqref{4-energy} and \eqref{4-StriQLW} hold for $u$ since they hold
uniformly for $u^k$.

\subsection{Reduction to Dispersive Estimate}\label{4-sec-RedDisp}

Here we show briefly how Proposition \ref{4-prop-localexst} follows
from Theorem \ref{4-thm-disp}.

Let us first recall some notions in \cite{SmTa05} which is
necessary for proceeding. Let $n=2$, and $\theta=\ep_0^{1/2}
\la^{-1/2}$ with $\la\gg 1$ stands for the frequency  and
$\ep_0\ll 1$ s.t. $\ep_0 \la\gg 1$, we use $\chi_{j,k,\om}$ to
denote the $L_x^\infty$-normalized wave packet supported in the
region (with $x_\om=x\cdot\om$ and $x'_\om$ be the given
orthonormal coordinates) \beeq
\label{4-wpsupp}T_{j,k,\om}=\{(x,t):\ |x_\om-t-k \la^{-1}|\le
\la^{-1},\ |x'_\om - j (\ep_0 \la)^{-1/2}|\le (\ep_0 \la)^{-1/2},\
|t|\le 2\}\eneq Note that for simplicity, we write here all the
quantities with respect to the flat metric, and this is sufficient
for us as explained at the beginning of Section
\ref{4-sec-sketch}. Precisely, \beeq
\label{4-wpdefi}\chi_{j,k,\om}=\la^{-1} T_\la(\delta(x_\om-t-k
\la^{-1}) W),\eneq where $T_\la$ is the convolution with a
spatially localized function $\psi_\la(x)=\la^n \psi(\la x)$, and
$W=W_0((\ep_0 \la)^{\frac{1}{2}}(x'_\om-j (\ep_0 \la)^{-1/2}))$.
The index $\om$, which stands for the initial orientation of the
wave packet at $t = -2$, varies over a maximal collection of
approximately $\theta^{-1}$ unit vectors separated by at least
$\theta$.

If we define the ``radial" wave packet \beeq
 \label{4-radialwp}\chi_{j,k}=\sum_{\om} \chi_{j,k,\om},\eneq
then as in \cite{SmTa05}, Proposition \ref{4-prop-localexst}
follows from the dispersive estimate in Theorem \ref{4-thm-disp}
for the superposition of radial wave packet.

We outline here how Theorem \ref{4-thm-disp} yields Proposition
\ref{4-prop-localexst}, for the details of the Propositions we used,
one should consult the content in \cite{SmTa05}. Firstly,
Proposition \ref{4-prop-localexst} is the consequence of the
following result which is similar to Proposition 7.2 in
\cite{SmTa05}. Let $S_\la$ or $S_{<\la}$ be the Littlewood-Paley
projector at or below the frequency $\la$, and $g_\la=S_{<\la}g$.
\begin{prop}\label{4-Prop-Pr7.2}
  Let $\ep_o \la\gg 1$, Then for
each $(u_0, u_1)\in H^1 \times L^2$, there exists a function $u_\la$
in $C^\infty([-2,2]\times \bR^2$ with
$$\mathrm{supp}\ (\widehat{u_\la(t, \cdot)}(\xi))\subset \{\xi: \la/8\le |\xi|\le 8 \la\}, $$
such that \beeq\label{4-eqn-approxslt}\| \Box_{g_\la} u_\la\|_{L^1_t
L_x^2}\lesssim \ep_0( \|u_0\|_{H^1}+\|u_1\|_{L^2})\ ,\eneq
\beeq\label{4-eqn-data} u_\la(-2)=S_\la u_0,\ \pt
 u_\la(-2)=S_\la u_1\ ,\eneq
 and such that the following Strichartz estimate holds for
 $r>\frac{1}{2}$
\beeq\label{4-eqn-StriSingFreq}\|S_\la u_\la\|_{L^2_t
L^\infty_x}\lesssim \ep_0^{-\frac{5}{4}} \la^{r-1}
(\|u_0\|_{H^1}+\|u_1\|_{L^2})\ .\eneq
\end{prop}

Now we use Theorem \ref{4-thm-disp} to give the proof of Proposition
\ref{4-Prop-Pr7.2}.
Let $u_{j,k,\om}=\theta^{\frac{1}{2}}\chi_{j,k,\om}$.
Then by Proposition 8.7 in \cite{SmTa05}, for any radial
$(u_0,u_1)\in H^1\times L^2$, there exists a function of form
 $$u=\sum_{j,k,\om} a_{j,k} u_{j,k,\om}$$ such that the equality
 \eqref{4-eqn-data} holds.
Moreover, by Proposition 8.4 in \cite{SmTa05}, we have \beeq
\label{4-prop8.7} \ep_0^{-1}\|\Box_{g_\la} S_\la u\|_{L^1_t L_x^2} +
\|\pa S_\la u\|_{L^\infty_t L^2_x}\lesssim
\|a_{j,k}\|_{l^2_{j,k,\om}}\lesssim \|u_0\|_{H^1}+\|u_1\|_{L^2}\
.\eneq
Now if we apply Theorem \ref{4-thm-disp} to $u$ here, we get that
\begin{eqnarray*}
\|S_\la u\|_{L^2_t L^\infty_x}&=&
\theta^{\frac{1}{2}}\|\sum_{j,k,\om} a_{j,k}
\chi_{j,k,\om}\|_{L^2_tL^\infty_x}\\
&=&\theta^{\frac{1}{2}}\|\sum_{j,k} a_{j,k}
\chi_{j,k}\|_{L^2_t L^\infty_x}\\
& \lesssim & \theta^{\frac{1}{2}} \times \ep_0^{-\frac{7}{4}}
(\ln \la)^{\frac{3}{2}}\|a_{j,k}\|_{l_{j,k}^2}\\
&\lesssim&\theta \times \ep_0^{-\frac{7}{4}} (\ln \la)^{\frac{3}{2}}
\|a_{j,k}\|_{l_{j,k,\om}^2}\\
&\lesssim& \ep_0^{-\frac{5}{4}} \la^{-\frac{1}{2}}(\ln
\la)^{\frac{3}{2}} (\|u_0\|_{H^1}+\|u_1\|_{L^2})\ ,
\end{eqnarray*} This is just the required Strichartz estimate at frequency
$\la$ \eqref{4-eqn-StriSingFreq}. Thus we complete the proof of
Proposition \ref{4-Prop-Pr7.2}.

\section{Dispersive Estimate}\label{4-sec-sketch}
In this section, we reduce Theorem \ref{4-thm-disp} to the proof of
Proposition \ref{4-PROP-A_1}, \ref{4-PROP-A_2} and \ref{4-PROP-over}
below, which deal with three sets of $(j,k)$ separately.

Based on the estimate of the Hamiltonian flow in \cite{SmTa05},
without loss of generality, we need only to give the proof of
Theorem \ref{4-thm-disp} for the flat metric.

In the process of the study, we find that one should deal with three
cases separately. Define the following subsets of the indices
$(j,k)$ in $\bZ^2$,
$$A_1=\{(j,k)|\ j^2 (\ep_0 \la)^{-1}+k^2 \la^{-2}\gg 1\}$$
$$A_2=\{(j,k)|\ j^2 (\ep_0 \la)^{-1}+k^2 \la^{-2}\ll 1,\ |j|\gg 1\}$$
$$A_3=\{(j,k)|\ j^2 (\ep_0 \la)^{-1}+k^2 \la^{-2}\ll 1,\ |j|\lesssim 1\}$$
We will prove Theorem \ref{4-thm-disp} for $(j,k)\in A_i$
separately.

In the case of $A_1$ and $A_2$, we have $\chi_{j,k}\lesssim 1$ in
principle which will be clear in Proposition \ref{4-PROP-over}, and
hence the dispersive estimate reduced to overlap estimate of the
wave packet as in Proposition 10.1 of \cite{SmTa05}.

Let $P_i=(t_i,x_i)$, and define $$N_i(P_1,P_2)=|\{(j,k)\in A_i \ |\
\chi_{j,k}(P_1) \chi_{j,k}(P_2)\neq 0 \}|$$ Then based on the
estimate Proposition 9.2 in \cite{SmTa05}, we can get the estimate
of $N_i(P_1,P_2)$ for $i=1,2$.
\begin{prop}\label{4-PROP-A_1}We have
\beeq \label{4-A1est} N_1(P_1,P_2)\lesssim \ep_0^{-3/2}
|t_1-t_2|^{-1/2}\eneq
\end{prop}

\begin{prop}\label{4-PROP-A_2}We have
\beeq \label{4-A2est}N_2(P_1,P_2)\lesssim \ep_0^{-1}
|t_1-t_2|^{-1}\eneq
\end{prop}

For the remained case $A_3$, the previous argument doesn't work.
Instead, we prove the $L^\infty$ estimate for the $\chi_{j,k}$ with
$(j,k)\in A_3$.
\begin{prop}\label{4-PROP-over}If $|j|\lesssim 1$, $|k|\lesssim \la$,
\beeq \label{4-overlap}|\chi_{j,k}(t,x)|\lesssim
\theta^{-1}\langle\la x\rangle^{-1/2}\langle\la |x|+1-|k+ \la
t|\rangle^{-1/2} H(\la |x|+1-|k+ \la t|)\eneq If  $|j|\lesssim 1$
and $\la \ll |k| \lesssim\theta^{-2}$, then
$|\chi_{j,k}(t,x)|\lesssim \ep_0^{-1}$. Else, 
$|\chi_{j,k}(t,x)|\lesssim 1$.
\end{prop}

\begin{coro}\label{4-coro}
$$\chi_{j,k}\lesssim \left\{\begin{array}{ll}\ep_0^{-1}&(j,k)\in
A_1\\1&(j,k)\in A_2\end{array}\right.$$
\end{coro}

By the previous result, we can prove Theorem \ref{4-thm-disp}
directly. In fact, by Proposition \ref{4-PROP-A_1} and Corollary
\ref{4-coro}, we have (as in Proposition 10.1 of \cite{SmTa05})
\beeq \label{4-A1final}\|\sum_{(j,k)\in A_1}
a_{j,k}\chi_{j,k}\|_{L^2_t L^\infty}\lesssim \ep_0^{-\frac{7}{4}}
(\ln\la)^{\frac{1}{2}}\|a_{j,k}\|_{l_{j,k}^2}.\eneq

And Proposition \ref{4-PROP-A_2} and Corollary \ref{4-coro} yields
\beeq \label{4-A2final}\|\sum_{(j,k)\in A_2}
a_{j,k}\chi_{j,k}\|_{L^2_t L^\infty}\lesssim \ep_0^{-\frac{1}{2}}
(\ln\la)^{\frac{3}{2}}\|a_{j,k}\|_{l_{j,k}^2}.\eneq

By Proposition \ref{4-PROP-over}, we have
\begin{prop}\label{4-A_3}
  \beeq \label{4-A3final}\|\sum_{(j,k)\in A_3} a_{j,k}\chi_{j,k}\|_{L^2_t
L^\infty}\lesssim \ep_0^{-\frac{1}{2}} \ln\la\
\|a_{j,k}\|_{l_{j,k}^2}.\eneq
\end{prop}

Thus, by \eqref{4-A1final}, \eqref{4-A2final}, \eqref{4-A3final},
Theorem \ref{4-thm-disp} is finally reduced to the proof of
Proposition \ref{4-PROP-A_1}, \ref{4-PROP-A_2} and
\ref{4-PROP-over}.

We give the proof of Proposition \ref{4-A_3} now.\\
\noindent {\bf Proof of Proposition \ref{4-A_3}:} Without loss of
generality, let $j=0$ and
$$f_k(t,m)=\langle m\rangle^{-\frac{1}{2}}\langle m+1-|k+t
\la|\rangle^{-\frac{1}{2}}H( m+1-|k+t \la|)\ ,$$ then by
\eqref{4-overlap},
\eqref{4-A3final} is reduced to the proof of \beeq
\label{4-reduceA3}\|\sum_{|k|\lesssim \la} a_k f_k\|_{L^2_t
L^\infty_m}\lesssim \la^{-\frac{1}{2}} \ln\la\ \|a_{k}\|_{l_{k}^2}\
. \eneq Since $f_k\ge 0$, we may assume $a_k\ge 0$ without loss of
generality. Let $$f_k^{\pm 1}(t,m)=f_k(t,m) H(\pm (k+t\la \mp 2))$$
and $f_k^0=f_k-f_k^1-f_k^{-1}$.

Since for any fixed $t$, there is finite $k$($|k+ t \la|<2$) such
that $f_k^0$ nonzero. The estimate for $f_k^0$ follows
directly($|f_k^0|\le 1$),
$$\|\sum_{|k|\lesssim \la} a_k f_k^0\|^2_{L^2_t L^\infty_m}\lesssim
\la^{-1} \sum a_{k}^2 .$$ Thus 
we need only to prove \eqref{4-reduceA3} for $f_k^1$ with $a_k\ge
0$, by symmetry.

Divide the time interval $[-2,2]$ into
$I_i=[\frac{i}{\la},\frac{i+1}{\la}]$ with $|i|\lesssim \la$. Then
for any $t \in I_i$, \begin{equation}\label{4-superposition}
\sum_{k} a_k f_k^1(t,m)\lesssim \sum_{1\le k+i\le m+1} a_k \langle
m\rangle^{-\frac{1}{2}}\langle m+1-k-i\rangle^{-\frac{1}{2}}=R(m).
\end{equation}  Let $m_i\lesssim \la$ be the point such that$$\|R(m)\|_{L^\infty_m}=R(m_i),$$
then
\begin{eqnarray*}
\|\sum_{k} a_k f_k^1(t,m)\|_{L^2_t L^\infty_m}^2& \lesssim & \sum_i
\left|\sum_{1\le k+i\le m_i+1} a_k \langle
m_i\rangle^{-\frac{1}{2}}\langle
m_i+1-k-i\rangle^{-\frac{1}{2}}\right|^2 \la^{-1}\\
&\lesssim &\la^{-1} \sum_i \langle m_i\rangle^{-1} (\sum_{k} a_k^2)
(\sum_{k} \langle m_i+1-k-i\rangle^{-1})\\
&\lesssim & \la^{-1} \ln\la \sum_{1\le k+i\le m_i+1} a_k^2
(m_i+1)^{-1}\\
&\lesssim & \la^{-1} \ln\la \sum_k \left(a_k^2 \sum_{i\ge
1-k}(k+i)^{-1}\right)\\
&\lesssim & \la^{-1} (\ln\la)^2 \sum_k a_k^2.
\end{eqnarray*}
This is just \eqref{4-reduceA3} for $f_k^1$. {\hfill  {\vrule
height6pt width6pt depth0pt}\medskip}

\section{Overlap estimates}\label{4-sec-overlap}
We first recall Proposition 9.2 in \cite{SmTa05} which is
essential for the proof of the overlap estimates. Let
$P_i=(t_i,x_i)$, $t_2>t_1$ and $t=t_2-t_1$. Note that
$$m=\max_{\om} ((x_2-x_1)\cdot \om - t)=|x_2-x_1|-t,$$ where the
maximum is attained at $\om=\alpha:=\frac{x_2-x_1}{|x_2-x_1|}$. We
define
$$N_{\la}(P_1,P_2) = | \{ (j,k,\om)\ |\
\chi_{j,k,\om}(P_1) \chi_{j,k,\om}(P_2)\neq 0\ \} |\ .$$

\begin{lem}[Proposition 9.2 in \cite{SmTa05}]\label{4-lem-Prop92}
\begin{equation}\label{4-Prop92}
  N_\la(P_1,P_2)\lesssim\left\{\begin{array}{ll}
    \theta^{-1}\langle\la  m\rangle^{-\frac{1}{2}}\langle\la
    t\rangle^{-\frac{1}{2}}&
    -4 \la^{-1}\le m\le \min(2t,c (\ep_0 \la)^{-1} t^{-1})\\
    \theta^{-1}\langle\la  m\rangle^{-1}&2t\le m\le c (\ep_0
    \la)^{-\frac{1}{2}}\\
    0& \mathrm{else}
  \end{array}\right.
\end{equation}
Precisely, let $t\ge \la^{-1}$ and
\begin{equation}\label{4-Ala-defi} A_\la=\{\om\in S^1 \ |\
|(x_2-x_1)\cdot \om-(t_2-t_1)|\lesssim \la^{-1}, \
|x_2-x_1-(t_2-t_1)\om |\lesssim (\ep_0 \la)^{-\frac{1}{2}}
\}\end{equation} then the estimate of $N_\la(P_1,P_2)$ follows from
the estimate of the area of $A_\la$ by the inequality
$$N_\la(P_1,P_2)\lesssim \theta^{-1}|A_\la(P_1,P_2)|\ .$$ Moreover,
if $|m|\lesssim \la^{-1}$ and $t \gtrsim \la^{-1}$,
\begin{equation}\label{4-Ala1} A_\la\subset \{|\om-\alpha|\lesssim
(\la t)^{-\frac{1}{2}}\}\ .
\end{equation}
If $\la^{-1}\lesssim m \lesssim \min(t, (\ep_0 \la)^{-1} t^{-1})$,
\begin{equation}\label{4-Ala2} A_\la\subset \{|\om-\alpha|\simeq
m^{\frac{1}{2}} t^{-\frac{1}{2}}\}\ .
\end{equation}
\end{lem}

Now we are ready to give the proof of Proposition \ref{4-PROP-A_1}
and
\ref{4-PROP-A_2}.\\
\noindent {\bf Proof of Proposition \ref{4-PROP-A_1}:} Without loss
of generality, we assume $P_i=(t_i,r_i,0)$ with $r_2\ge r_1\gg 1$
and $t=t_2-t_1\ge 0$. Denote the spatial clockwise rotation of $P_2$
with angle $\om$ by $P_2^{\om}$. Define for $k\theta\in [0,\pi]$
$$m_k:=m(P_1,P_2^{k\theta})=\sqrt{(r_2-r_1)^2+2 r_1 r_2 (1-\cos(k\theta))}-t\ge r_2-r_1-t\ .$$

First, for $t\gg (\ep_0 \la)^{-\frac{1}{2}}$, we have
$N_\la(P_1,P_2^{k\theta})$ nonzero only if $m_k\in(-4\la^{-1}, c
(\ep_0 \la t)^{-1})$ by \eqref{4-Prop92}. Thus w.l.o.g, we may
assume $r_2-r_1-t\le c (\ep_0 \la t)^{-1}$.
We consider separately the cases $|r_2-r_1-t|\le c (\ep_0 \la
t)^{-1}$ and $r_2-r_1-t\le - c (\ep_0 \la t)^{-1}$.

For the first case, $|r_2-r_1-t|\le c (\ep_0 \la t)^{-1}$, we have
from $m_k\lesssim (\ep_0 \la t)^{-1}$
$$1-\cos(k\theta)\le \frac{(t+ c (\ep_0 \la t)^{-1})^2-(r_2-r_1)^2}{2 r_1 r_2}\lesssim
\frac{t\cdot (\ep_0 \la t)^{-1}}{r_1 r_2}\ll (\ep_0 \la )^{-1}$$
thus $k\theta\ll (\ep_0 \la )^{-\frac{1}{2}}$ and $k\ll \ep_0^{-1}$.
So $$N_1(P_1,P_2)\lesssim \ep_0^{-1} \max_k
\{N_\la(P_1,P_2^{k\theta}\}\lesssim \ep_0^{-1} \cdot \theta^{-1}
\langle\la
    t\rangle^{-\frac{1}{2}}\lesssim \ep_0^{-\frac{3}{2}} t^{-\frac{1}{2}}$$

For the second case $r_2-r_1-t\le -c (\ep_0 \la t)^{-1}$, let $k_1$
$k_2$ be the number s.t. $m_{k_1}=c (\ep_0 \la t)^{-1}$ and
$m_{k_2}=-4 \la^{-1}$. Then
\begin{eqnarray*}
1-\cos(k_1 \theta)&=& \frac{(t+ c (\ep_0 \la
t)^{-1})^2-(r_2-r_1)^2}{2 r_1 r_2}\\
&\simeq& \frac{t(t+r_1-r_2+c(\ep_0 \la t)^{-1})}{r_1 r_2} \\
&\gtrsim& (\ep_0 \la r_1 r_2)^{-1}.\end{eqnarray*} On the other
hand,
$$1-\cos(k_1 \theta)= \frac{(t+ c (\ep_0 \la
t)^{-1})^2-(r_2-r_1)^2}{2 r_1 r_2}\lesssim \frac{t^2}{r_1 r_2}\ll
1$$ and thus $1\gg k_1 \theta \gtrsim (\ep_0 \la r_1
r_2)^{-\frac{1}{2}}$. Note that
$$\cos(k_2 \theta)-\cos(k_1 \theta)=\frac{(t+ c (\ep_0 \la t)^{-1})^2-(t-4\la^{-1})^2}{2 r_1 r_2}\simeq (\ep_0 \la r_1 r_2)^{-1}$$
and $$\cos(k_2 \theta)-\cos(k_1 \theta)=\int_{k_2 \theta}^{k_1
\theta} \sin x d x\simeq (k_1-k_2)k_1 \theta^2\ .$$ So we have
$$
k_1-k_2\simeq (\ep_0 \la r_1 r_2)^{-1} \theta^{-2} k_1^{-1}\lesssim
(\ep_0 \la r_1 r_2)^{-\frac{1}{2}} \theta^{-1}\ll \ep_0^{-1}
$$ and hence
$N_1(P_1,P_2)\lesssim \ep_0^{-\frac{3}{2}} t^{-\frac{1}{2}}$ as
before.

Secondly, for $t\lesssim (\ep_0 \la)^{-\frac{1}{2}}$, we have
$N_\la(P_1,P_2^{k\theta})$ nonzero only if $m_k\in(-4\la^{-1}, c
(\ep_0 \la )^{-\frac{1}{2}})$ by \eqref{4-Prop92}. Then from
$m_{k}\le c (\ep_0 \la )^{-\frac{1}{2}}$,
$$1-\cos(k\theta)= \frac{(t+ c (\ep_0 \la )^{-\frac{1}{2}})^2-(r_2-r_1)^2}{2 r_1
r_2}\lesssim  \frac{(\ep_0 \la )^{-1}}{r_1 r_2} \ll (\ep_0 \la
)^{-1}$$ So $k\ll \theta^{-1} (\ep_0 \la )^{-\frac{1}{2}}=
\ep_0^{-1}$, and hence $N_1(P_1,P_2)\lesssim \ep_0^{-\frac{3}{2}}
t^{-\frac{1}{2}}$ as before. This completes the proof of Proposition
\ref{4-PROP-A_1}.{\hfill  {\vrule height6pt width6pt
depth0pt}\medskip}

\noindent {\bf Proof of Proposition \ref{4-PROP-A_2}:} We use the
same notation as in the proof of Proposition \ref{4-PROP-A_1}. Note
that since $P_i\in A_2$, thus we may assume $1\gtrsim r_2\ge r_1 \gg
(\ep_0 \la)^{-\frac{1}{2}}$ and $t\ge 0$. Let $\alpha_2=k
\theta\in[0,\pi]$ and note that
$P_2^{\alpha_2}-P_1=(t,r,0)^{\alpha}$ with
\begin{equation}\label{4-rexpress}
r^2=(r_2-r_1)^2+2 r_1 r_2 (1-\cos(k\theta))\end{equation} and
\begin{equation}\label{4-alphaexpress}
r \sin \al=r_2 \sin \al_2\ .\end{equation}

Since $$N_2(P_1,P_2)
\lesssim |\{(j,k,\om)\ |\ \chi_{j,k,\om}(P_1) \neq 0 \}|\lesssim
\theta^{-1}\ ,$$ we may assume $t\gg (\ep_0 \la)^{-\frac{1}{2}}$
w.l.o.g.. Since we are restricted in $A_2$, where $|j|\gg 1$, we can
modify the definition of $A_\la$ as follows, $$A_\la=\{\om\in S^1 \
|\  |\sin\om|\gg (\ep_0 \la)^{-\frac{1}{2}} r_1^{-1},\ \exists j,k\
\mathrm{s.t.}\ \chi_{j,k,\om}(P_1)\chi_{j,k,\om}(P_2)\neq 0\}$$

Note that $m=r-t$ depend on $\al_2$, \beeq
\label{4-mvaries}\partial_{\al_2} m=\frac{r_1 r_2 \sin
\al_2}{r}=r_1 \sin \al.\eneq

If $(\ep_0 \la t)^{-\frac{1}{2}} t^{-\frac{1}{2}}\lesssim (\ep_0
\la )^{-\frac{1}{2}} r_1^{-1}$, i.e., $r_1\lesssim t$, then by
\eqref{4-Ala1} and \eqref{4-Ala2}, for any $\om\in A_\la$, one has
$|\om-\al|\lesssim (\ep_0 \la )^{-\frac{1}{2}} r_1^{-1}$. Thus we
have $|\sin \al |\gtrsim (\ep_0 \la )^{-\frac{1}{2}} r_1^{-1}$ in
case of $A_\la$ nonempty. So
\begin{eqnarray*}\sum_{\al_2}|A_{\la}(P_1,P_2^{\al_2})|&\lesssim&
(\la t)^{-\frac{1}{2}} \frac{c (\ep_0 \la t)^{-1}-(-4
\la^{-1})}{\sup_{\al_2}|\partial_{\al_2}m|}\\
&\lesssim& (\la t)^{-\frac{1}{2}} \frac{(\ep_0 \la t)^{-1}}{(\ep_0
\la)^{-\frac{1}{2}}}\\
&\lesssim & (\la t)^{-\frac{1}{2}}\end{eqnarray*} and
$$N_2(P_1,P_2)\lesssim \theta^{-1} (\la
t)^{-\frac{1}{2}} =(\ep_0 t)^{-\frac{1}{2}}\lesssim(\ep_0 t)^{-1}\
.$$

Else if $r_1\gg t$, we assume $r_2-r_1-t\le c (\ep_0 \la t)^{-1}$
w.l.o.g.. so $r_2\simeq r_1$. Let $k_1$ $k_2$ be the number s.t.
$m_{k_1}=c (\ep_0 \la t)^{-1}$ and $m_{k_2}=\max(-4 \la^{-1},
r_2-r_1-t)$. We claim that \beeq \label{4-claim}k_1-k_2\lesssim
(\ep_0 r_1)^{-1}\ .\eneq 
We consider separately the cases $|r_2-r_1-t|\le c (\ep_0 \la
t)^{-1}$ and $r_2-r_1-t\le - c (\ep_0 \la t)^{-1}$.

For the first case, $|r_2-r_1-t|\le c (\ep_0 \la t)^{-1}$, we have
\begin{eqnarray*}1-\cos(k_1 \theta)&\le& \frac{(t+ c (\ep_0 \la
t)^{-1})^2-(r_2-r_1)^2}{2 r_1 r_2}\\
&\lesssim& \frac{t\cdot (\ep_0 \la t)^{-1}}{r_1 r_2}\\
&\lesssim& (\ep_0 \la )^{-1} r_1^{-2}\ .\end{eqnarray*} thus
$k_1\theta\lesssim (\ep_0 \la )^{-\frac{1}{2}}r_1^{-1}$.

For the second case $r_2-r_1-t\le -c (\ep_0 \la t)^{-1}$,
\begin{eqnarray*}1-\cos(k_1 \theta)&=& \frac{(t+ c (\ep_0 \la
t)^{-1})^2-(r_2-r_1)^2}{2 r_1 r_2}\\
& \simeq & \frac{t(t+r_1-r_2+c(\ep_0 \la t)^{-1})}{r_1 r_2}
\\
&\gtrsim & (\ep_0 \la )^{-1} r_1^{-2}.\end{eqnarray*} On the other
hand,
$$1-\cos(k_1 \theta)= \frac{(t+ c (\ep_0 \la
t)^{-1})^2-(r_2-r_1)^2}{2 r_1 r_2}\lesssim \frac{t^2}{r_1^2}\ll 1$$
and thus $1\gg k_1 \theta \gtrsim (\ep_0 \la
)^{-\frac{1}{2}}r_1^{-1}$. Note that
$$\cos(k_2 \theta)-\cos(k_1 \theta)=\frac{(t+ c (\ep_0 \la t)^{-1})^2-(t-4\la^{-1})^2}{2 r_1 r_2}\simeq (\ep_0 \la )^{-1} r_1^{-2},$$
we have $(k_1-k_2) \theta \lesssim (\ep_0 \la
)^{-\frac{1}{2}}r_1^{-1}$. This proves the claim \eqref{4-claim}.

Now we are ready to estimate $N_2(P_1,P_2)$. If $(\ep_0
\la)^{-\frac{1}{2}} r_1^{-1}\ll (\la t)^{-\frac{1}{2}}$, then
\begin{eqnarray*}N_2(P_1,P_2)&\lesssim& \theta^{-1}
\sum_{\al_2}|A_{\la}(P_1,P_2^{\al_2})|\\
&\lesssim& \theta^{-1}(k_1-k_2)(\la t)^{-\frac{1}{2}}\\
&\lesssim& \theta^{-2}(\la t)^{-1}=(\ep_0 t)^{-1}\ .
\end{eqnarray*}
Else if $(\ep_0 \la)^{-\frac{1}{2}} r_1^{-1}\gtrsim (\la
t)^{-\frac{1}{2}}$, let $m_0$ be s.t.
$m_0^{\frac{1}{2}}t^{-\frac{1}{2}}=(\ep_0 \la)^{-\frac{1}{2}}
r_1^{-1}$ and $k_0$ s.t. $m_{k_0}=\max(m_0,r_2-r_1-t)$. then for
$k\in [k_2,k_0]$, $$\partial_{\al_2}m=r_1 \sin \al\gtrsim (\ep_0
\la)^{-\frac{1}{2}}$$
\begin{eqnarray*}\sum_{\al_2}|A_{\la}(P_1,P_2^{\al_2})|
&\lesssim& \frac{m_0}{(\ep_0 \la)^{-\frac{1}{2}}}(\la
t)^{-\frac{1}{2}}\\
&\lesssim &(\la t)^{-\frac{1}{2}}\\
&\lesssim & (\ep_0 \la)^{-\frac{1}{2}} r_1^{-1} \lesssim (\ep_0
\la)^{-\frac{1}{2}} t^{-1}.\end{eqnarray*}
 For $k\in
[k_0,k_1]$,
\begin{eqnarray*}
\sum_{\al_2}|A_{\la}(P_1,P_2^{\al_2})| &\lesssim& (k_1-k_0)( \la
m_0)^{-\frac{1}{2}}(\la t)^{-\frac{1}{2}}\\
&\lesssim& \theta^{-1} (\ep_0 \la )^{-\frac{1}{2}}r_1^{-1} \la^{-1}
(m_0 t)^{-\frac{1}{2}}\\
&=&(\theta \la t)^{-1}= (\ep_0 \la)^{-\frac{1}{2}}
t^{-1}.\end{eqnarray*} Thus $$N_2(P_1,P_2)\lesssim
\theta^{-1}(\ep_0 \la)^{-\frac{1}{2}} t^{-1}=(\ep_0 t)^{-1}.\ $$
This completes the proof of Proposition \ref{4-PROP-A_2}.{\hfill
{\vrule height6pt width6pt depth0pt}\medskip}

\section{$L^\infty$ estimate for radial wave packet}\label{4-sec-est-ra-wp}
Since $\chi_{j,k,\om}(t,x)=\chi_{j,k-t\la,\om}(0,x)$, we may assume
$t=0$ and $j,k\ge 0$ w.l.o.g..

We first show that \beeq \label{4-disjoint}T_{j,k,0}\cap
T_{j,k,\theta}=\emptyset\ \mathrm{for}\ j\gg 1 \ \mathrm{or}\ k\gg
\theta^{-2}.\eneq If $j=0$, $k\gg 1$,  and $T_{j,k,0}\cap
T_{j,k,\theta}\neq\emptyset$, then \begin{equation}\label{4-klarge}
(k-1)\la^{-1}\tan \frac{\theta}{2}\le \frac{1}{2}(\ep_0
\la)^{-\frac{1}{2}}\ ,\end{equation} i.e., $k\lesssim \theta^{-2}$.
If $k=0$, $j\gg 1$,  and $T_{j,k,0}\cap
T_{j,k,\theta}\neq\emptyset$, then
$$(j-1)(\ep_0 \la)^{-\frac{1}{2}} \tan \frac{\theta}{2}\le
\frac{1}{2}(\la)^{-1}\ ,$$ i.e., $j\lesssim 1$.

If $j,k\ge 1$, and $T_{j,k,0}\cap T_{j,k,\theta}\neq\emptyset$, let
$(x,y):=((j-1)(\ep_0 \la)^{-\frac{1}{2}}, (k+1)\la^{-1})\in
T_{j,k,0}$, then $(x,y)^{\theta}\in T_{j,k,0}$. Thus
$$\left\{\begin{array}{lcl}
y\cos \theta-x \sin \theta&\ge& (k-1)\la^{-1}\\
x\cos \theta+y \sin \theta&\le& (j+1)(\ep_0 \la)^{-\frac{1}{2}}.
\end{array}\right.$$
From the first inequality, we have $$(j-1)(\ep_0 \la)^{-\frac{1}{2}}
\sin \theta\lesssim \la^{-1}\ ,$$ so $j\lesssim 1$. Then by the
second inequality, $$(k+1)\la^{-1}\sin \theta\lesssim(\ep_0
\la)^{-\frac{1}{2}}\ ,
$$ thus $k\lesssim \theta^{-2}$. Combined these observations, we get
\eqref{4-disjoint}.

Similar argument will yield \beeq \label{4-disjoint2}T_{j,k,0}\cap
T_{j,k, M\theta}=\emptyset\ \mathrm{for}\ j\gg 1 \ \mathrm{or}\ k\gg
M^{-1}\theta^{-2}\ ,\eneq for $M$ s.t. $M\theta\ll 1$. Then we have
$|\chi_{j,k}(0,x)|\lesssim 1$ for $|j|\gg 1$ or $|k|\gg
\theta^{-2}$, and $|\chi_{j,k}(0,x)|\lesssim \ep_0^{-1}$ for
$|j|\lesssim 1$ and $\la \lesssim |k|\lesssim \theta^{-2}$.

It remains to consider the case $j\lesssim 1$ and $k\lesssim \la$
now. Note that the estimate for $j\lesssim 1$ can be reduced to
the counterpart for $j=0$, Proposition \ref{4-PROP-over} follows
from the following Lemma \ref{4-lemover}

\begin{lem}\label{4-lemover}
Let $0\le k\le \la$, we have \beeq
\label{4-overlap1}|\chi_{0,k}(0,x)|\lesssim \theta^{-1}\langle\la
x\rangle^{-1/2}\langle\la |x|+1-k\rangle^{-1/2} H(\la
|x|+1-k).\eneq
\end{lem}
\begin{prf}
If $|x|< (k-1)\la^{-1}$, then $\chi_{0,k,\om}(0,x)=0$ for any $\om$
and hence $\chi_{0,k}(0,x)=0$. Thus for the proof of
\eqref{4-overlap1}, we need only to show for $|x|\ge (k-1)\la^{-1}$,
\beeq \label{4-overlap2}|\chi_{j,k}(0,x)|\lesssim
\theta^{-1}\langle\la x\rangle^{-1/2}\langle\la
|x|+1-k\rangle^{-1/2}\ . \eneq

If $k\lesssim 1$. For the case $|x|\lesssim\la^{-1}$, we use the
trivial bound $|\chi_{0,k}(0,x)|\lesssim \theta^{-1}$. For $|x|\gg
(\ep_0 \la)^{-\frac{1}{2}}$, it is obviously that
$\chi_{0,k,\om}(0,x)=0$ for any $\om$. For the remained case
$\la^{-1}\ll |x|\lesssim (\ep_0 \la)^{-\frac{1}{2}}$, it's only need
to calculate the number of $l$ s.t. $(0,x)\in T_{0,k,l\om}$, denoted
by $\tau$. Then we have $$|x| \sin \tau \theta \lesssim \la^{-1}\
.$$ Thus $\sin \tau \theta \ll 1$ and hence $|x| \tau \theta
\lesssim \la^{-1}$, i.e., $\tau\lesssim \theta^{-1}(|x|\la)^{-1}$.

We consider for $k\gg 1$ now. Let $A(x)=\{\om |\ (0,x)\in
T_{0,k,\om}\}$. Since $(0,x)\in T_{0,k,\om}$, then $$|x\cdot \om-k
\la^{-1}|\lesssim \la^{-1}, \ |x-k \la^{-1}\om |\lesssim (\ep_0
\la)^{-\frac{1}{2}}\ .$$ Compare it with the definition
\eqref{4-Ala-defi} of $A_\la(P_1,P_2)$, we get that $$A(x) \subset
A_\la((0,x),(-k \la^{-1}, 0))\ .$$ Thus
\begin{eqnarray*}
|\chi_{0,k}(0,x)| &\le& \sum_{|j \theta|\le \pi} |\chi_{0,k,
j\theta}(0,x)|\\
&\le& |\{j |\ j\theta\in [-\pi, \pi],\ \chi_{0,k, j\theta}(0,x)\neq
0\}|\\
&\lesssim& \theta^{-1}|A(x)|\\
&\lesssim& \theta^{-1}|A_\la(P_1,P_2)|\ ,
\end{eqnarray*}
where $P_2=(0,x)$ and $P_1=(-k \la^{-1}, 0)$.

By the notation at the beginning of Section \ref{4-sec-overlap}, we
have $t=k \la^{-1}\le 1$ and $m = |x| -k \la^{-1}$. Then by Lemma
\ref{4-lem-Prop92}, if $t\lesssim (\ep_0 \la)^{-\frac{1}{2}}$, i.e.,
$k\lesssim \theta^{-1}$, we have $$A_\la\lesssim
\left\{\begin{array}{ll} \langle\la
m\rangle^{-\frac{1}{2}}\langle\la t\rangle^{-\frac{1}{2}}=\langle\la
|x|-k\rangle^{-\frac{1}{2}}\langle k\rangle^{-\frac{1}{2}}&
    m\lesssim t\\
    \langle\la  m\rangle^{-1}=\langle\la
|x|-k\rangle^{-1}& t\lesssim m\lesssim (\ep_0
    \la)^{-\frac{1}{2}}.
\end{array}
\right.$$ In both cases, we have $$A_\la\lesssim \langle\la
|x|-k\rangle^{-\frac{1}{2}}\langle \la |x|\rangle^{-\frac{1}{2}}\
.$$ If $ (\ep_0 \la)^{-\frac{1}{2}}\ll t\le 1$, i.e.,
$\theta^{-1}\ll k\le \la$, we have that for $m\lesssim (\ep_0
\la)^{-1}t^{-1}\ll t$ (thus $\la |x|\lesssim k$)
$$A_\la\lesssim
\langle\la m\rangle^{-\frac{1}{2}}\langle\la
t\rangle^{-\frac{1}{2}}=\langle\la
|x|-k\rangle^{-\frac{1}{2}}\langle k\rangle^{-\frac{1}{2}}\lesssim
\langle\la |x|-k\rangle^{-\frac{1}{2}}\langle \la
|x|\rangle^{-\frac{1}{2}}\ .
$$

\end{prf}



\end{document}